\documentclass[a4paper,12pt]{article}
\usepackage{amsmath, amsfonts, amssymb}
\usepackage[english]{babel}
\usepackage {color}

\definecolor{mycolor1}{rgb}{0.020,0.388,0.757}
\definecolor{mycolor2}{rgb}{0.208,0.208,0.208}

\definecolor{mycolor1}{rgb}{0.020,0.388,0.757}
\definecolor{mycolor2}{rgb}{0.208,0.208,0.208}

\begin{document}
\begin{center}
Lauricella hypergeometric function and its application to the solution of the Neumann problem for a multidimensional elliptic equation with several singular coefficients in an infinite domain
\end{center}

\begin{center}
T.G.~Ergashev, Z.R.~Tulakova
\end{center}

\textbf{Abstract.} At present, the fundamental solutions of the multidimensional elliptic equation with the several singular coefficients are known and they are  expressed in terms of the Lauricella hypergeometric function of many variables. In this paper we study the Neumann problem for a
multidimensional elliptic equation with several singular coefficients in the infinite domain. Using the method of the integral energy the uniqueness of solution has been proved. In the course of proving the existence of the explicit solution of the Neumann problem, a differentiation formula, some adjacent and limiting relations for the Lauricella hypergeometric functions and the values of some multidimensional improper integrals are used.



\bigskip

\section{Introduction}

It is known that the theory of boundary value problems for degenerate equations and equations with singular coefficients is one of the central sections of the modern theory of partial differential equations, which is encountered in solving many important problems of applied nature, for example, \cite{{bers},{A2}}. Omitting the huge bibliography in which various local and nonlocal boundary value problems for equations of mixed type containing elliptic equations with singular coefficients are studied, we note the works that are most closely related to this work.

Let $\mathbb{R}_m$ be the $ m $ -dimensional Euclidean space
\, $ (m \geq 2) $, $ x: = \left(x_1, ..., x_m \right) $ - arbitrary point
in it and $ n $ is a natural number, and $ n \leq m $.
The $ \, 2^n $ -th part of the Euclidean space $\mathbb{R}_m$ is defined as follows:

\begin{equation} \label{eq4}
\Omega \equiv \Omega_m^{n+}= \left\{x \in \mathbb{R}_m: x_i>0,\,
i=1,...,n; \,-\infty<x_j<+\infty,\, j=n+1,...,m\right\}.
\end{equation}

Fundamental solutions have an essential role in studying partial differential equations. The explicit form of the fundamental solution makes it possible to correctly formulate the problem statement and to study in detail the various properties of the solution of the equation under consideration.
Fundamental solutions  of singular elliptic equations are directly connected with multiple hypergeometric functions, the number of variables of which is determined by the number of singular coefficients. Indeed, all fundamental solutions of the following elliptic  equation with $n$ singular coefficients
\begin{equation}
\label{eq3} E_{\alpha} ^{(m,n)} (u) \equiv
\sum\limits_{i=1}^{m}{\frac{{{\partial }^{2}}u}{\partial x_{i}^{2}}}+\sum\limits_{j=1}^{n}{\frac{2{{\alpha }_{j}}}{{{x}_{j}}}\frac{\partial u}{\partial {{x}_{j}}}}=0
\end{equation}
in the hyperoctant $\Omega$ are expressed \cite{ergsib} by the Lauricella hypergeometric function $F_A^{(n)}$ in $n$ variables \cite{appell}  where $m\ge 2$ is a dimension of the Euclidean space; $n\ge 1$ is a number of the singular coefficients;  $m\ge n;$ ${{\alpha}_{j}}$  are real constants and $0<2{{\alpha}_{j}}<1$ $(j=\overline{1,n})$.

In a recent paper \cite{DU}, the generalized Holmgren problem  for equation (\ref{eq3}) in some part of the first hyperoctant of the ball (in the finite domain) is written out in an explicit form and in case $m=3$ and $n=1$ the potential theory  in the domain bounded in a half-space is constructed \cite{erglob1}.

Relatively few works are devoted to the study of boundary value problems for spatial singular elliptic equations in infinite domains. We note the joint work of M.Salakhitdinov and A.Hasanov \cite{salhas}, in which solutions of the  Dirichlet and Holmgren problems for a multidimensional elliptic equation with one singular coefficient in the half-space were found in an explicit forms. For three- and four-dimensional elliptic equations with two \cite{kar1}, three \cite{kar2} and four \cite{ryskan} singular coefficients boundary value problems in infinite domains are investigated. Recently \cite{ergtulakova2021}, a solution to the external Dirichlet problem for an equation (\ref{eq3}) in the hyperoctant $\Omega$ is found explicitly.

In this paper, we study the Neumann problem for equation (\ref{eq3}) in the infinite domain $ \Omega $.  The outer Neumann problem for equation (\ref{eq3}) on the plane has been studied in details by many authors (see, for example, \cite{marichev}), therefore, for definiteness, we set $m>2$ in this paper.
Initially we introduce some formulas and definitions, then we proceed to solve boundary-value problem.

\section{Preliminaries}

Below we give some formulas for Euler gamma-function, Gauss hypergeometric function, Lauricella hypergeometric functions of three and more variables, which will be used in the next sections.

Let be ${N}$ set of the natural numbers : $N=\{1,2,3,...\}.$

It is known that the Euler gamma-function $\Gamma(a)$ has property \cite[ eq. 1.2(2)]{Erd}
$$\Gamma(a+m)=\Gamma(a) (a)_m.$$
Here $(a)_m$ is a Pochhammer symbol, for which the equality
$(a)_{m+n}=(a)_m(a+m)_n $ and its particular case $(a)_{2m}=(a)_m(a+m)_m$ is true.

For the Gamma function $\Gamma(a)$ the Legendre's duplication formula
\cite[ eq. 1.2(15)]{Erd} :
\begin{equation} \label{legendre}
\Gamma \left( {2a} \right) = {\frac{{2^{2a - 1}}}{{\sqrt {\pi}} } }\Gamma
\left( {a} \right)\Gamma \left( {a + {\frac{{1}}{{2}}}} \right)
\end{equation}
is valid.

A function
\begin{equation*}
F {\left(a,b;c;x\right)} ={\sum\limits_{k =
0}^{\infty} \frac{(a)_k(b)_k}{(c)_kk!}x^k }, \,c\neq 0,-1,-2,...
\end{equation*}
is known as the  Gaussian hypergeometric function.

Multiple Lauricella  hypergeometric function $F_A^{(n)}$ in $n\in N$ (real or complex) variables is defined as following \cite[eq. 1.4(1)]{srikarl}:
\[
F^{(n)}_A \left( {{
 {a,{\bf{b}};}
 {{\bf{c}};}
} {\bf{x}}}  \right)
=\sum\limits_{|{\bf{k}}|=0}^\infty\frac{(a)_{|\bf{k}|}(b_1)_{k_1}...(b_n)_{k_n}}{k_1!...k_n!(c_1)_{k_1}...(c_n)_{k_n}}x_1^{k_1}...x_n^{k_n}\]
\[\,\,\,\,\,\,\,\,\,\,\,\,\,\,\,\,\,\,\,\,\,\,\,\,\,\left[c_{i} \ne 0,-1,-2,...;\,i = \overline {1,n} ; \,\,|x_1|+...+|x_n|<1\right],\]
where

\[{\bf{b}}:=\left(b_1,...,b_n\right),\, {\bf{c}}:=\left(c_1,...,c_n\right);\]
\[{\bf{x}}:=\left(x_1,...,x_n\right);\,
{|\bf{k}|}:=k_1+...+k_n,\,\,k_1\geq 0, ..., k_n\geq 0.\]

The function $ F_{A}^{\left({n} \right)} $ satisfies the following differentiation formula \cite{appell}

\begin{equation} \label{differential}
{\frac{{\partial}} {{\partial x_{k}}} }F_{A}^{(n)} \left( {{
 {a,\,{\textbf{b}};}\,
 {{\textbf{c}};}
} \,{\textbf{x}}}  \right)
 = {\frac{{ab_{k}}} {{c_{k}}} }F_{A}^{(n)}  \left( {{
 {a+1,\,{\textbf{{b}}}_k;}\,
 {{\textbf{{c}}}_k;}
}\, {\textbf{x}}}  \right)
\end{equation}
and adjacent relation \cite{DU}
\begin{equation}
\label{adjacent}
{\sum\limits_{k = 1}^{n} {}} {\frac{{a_{k}}} {{b_{k}}} }x_{k} F_{A}^{(n)}
\left(
 a + 1,  \textbf{b}_{k}; \textbf{c}_{k}; \textbf{x}  \right)  = F_{A}^{(n)} {\left( {{
 {a + 1,\textbf{b};}\textbf{c};} \textbf{x}}  \right)} -F_{A}^{(n)} {\left( {{
 {a,\textbf{b};}\textbf{c};} \textbf{x}}  \right)},
\end{equation}
where ${\textbf{b}_k}$  and   ${\textbf{c}_k}$  are vectors obtained, respectively, from vectors $ {\textbf{b}} $ and $ {\textbf{c}} $ by increasing the $ k $-th component by one ($ k = \overline{1, n} $).

\textbf{Lemma 1}.\cite{{ergtulakova2021}, {montes}}\label{L3} Let $a, b_{k}$ and  $c_{k}$ are real numbers with $ c_k \ne 0,\, - 1,\, - 2,...$, $a > b_{1} + ... +
b_{n}$ and $c_k>b_k$ \, $(k=\overline{1,n})$. Then for $n = 1,2,...$ the following limiting relation holds true
\[{\mathop {\lim} \limits_{\varepsilon \to 0}}\varepsilon^{-b_1-...-b_n}
F_{A}^{(n)} \left( {a,{\bf{b}}; {\bf{c}}
;  1-\frac{z_1(\varepsilon)}{\varepsilon}, ...,1-\frac{z_n(\varepsilon)}{\varepsilon}}  \right)
\]
\begin{equation}
\label{eq16}
= \frac{1}{\Gamma(a)}\Gamma \left( {a - {\sum\limits_{k = 1}^{n} {b_{k}}} }
\right)  \prod_{k=1}^n\frac{\Gamma\left(c_k\right)}{\left[z_k(0)\right]^{b_k}\Gamma\left(c_k-b_k\right)},
\end{equation}
where $z_k(\varepsilon)$ are arbitrary functions with $z_k(0)\neq 0$.

\textbf{Lemma 2} \label{grad}  If $p_k,\, q_k,\, r_k,\, s,\,t,\, $ are real numbers and
\[
{p_{k} > 0,\,\,q_{k} > 0,\,\,r_{k} > 0,\,\,s > 0,\,\,\,0 <
{\frac{{p_{1}}} {{q_{1}}} } + ... + {\frac{{p_{n}}} {{q_{n}}} } - t < s\,}, \, k=\overline{1,n},
\]
then the following equality holds true
\[
{\underbrace {{\int\limits_{0}^{ + \infty}  {...{\int\limits_{0}^{ + \infty
} {}}} }} _{n}}{\frac{{x_{1}^{p_{1} - 1} ...x_{n}^{p_{n} - 1} dx_{1}
...dx_{n}}} {{\left[ {\left( {r_{1} x_{1}}  \right)^{q_{1}}  + ... + \left(
{r_{n} x_{n}}  \right)^{q_{n}}}  \right]^{t}\left[ {1 + \left( {r_{1} x_{1}
} \right)^{q_{1}}  + ... + \left( {r_{n} x_{n}}  \right)^{q_{n}}}
\right]^{s}}}}
\]

\begin{equation}
\label{eq17}
 = {\frac{{\Gamma \left( {{\displaystyle\frac{{p_{1}}} {{q_{1}}} }} \right)...\Gamma
\left( {{\displaystyle\frac{{p_{n}}} {{q_{n}}} }} \right)\Gamma \left( {{\displaystyle\frac{{p_{1}
}}{{q_{1}}} } + ... + {\displaystyle\frac{{p_{n}}} {{q_{n}}} } - t} \right)\Gamma \left(
{s + t - {\displaystyle\frac{{p_{1}}} {{q_{1}}} } - ... - {\displaystyle\frac{{p_{n}}} {{q_{n}}} }}
\right)}}{{q_{1} q_{2} ...q_{n} r_{1}^{p_{1} q_{1}}  ...r_{n}^{p_{n} q_{n}}
\Gamma \left( {{\displaystyle\frac{{p_{1}}} {{q_{1}}} } + ... + {\displaystyle\frac{{p_{n}}} {{q_{n}
}}}} \right)\Gamma \left( {s} \right)}}}.
\end{equation}

\textbf{Proof}.
Into (\ref{eq17}) we make a replacement of variables
\[
\left( {r_{1} x_{1}}  \right)^{q_{1} / 2} = r\cos \phi _{1} ,\,\,\,\,\,\,\,\,\,\,
\]
\[
\left( {r_{2} x_{2}}  \right)^{q_{2} / 2} = r\sin \phi _{1} \cos \phi _{2}
,
\]
\[
\left( {r_{3} x_{3}}  \right)^{q_{3} / 2} = r\sin \phi _{1} \sin \phi _{2}
\cos \phi _{3} ,
\]
\[
................................
\]
\[
\left( {r_{n - 1} x_{n - 1}}  \right)^{q_{n - 1} / 2} = r\sin \phi _{1} \sin
\phi _{2} ...\sin \phi _{n - 2} \cos \phi _{n - 1} ,
\]
\[
\left( {r_{n} x_{n}}  \right)^{q_{n} / 2} = r\sin \phi _{1} \sin \phi _{2}
...\sin \phi _{n - 2} \sin \phi _{n - 1}
\quad
{\left[ {r \ge 0,\,\,0 \le \phi _{k} \le \pi / 2,\,\,k = \overline {1,n} \,}
\right]}
\]
and using the values of the following famous integrals \cite[eq. 1.5(16) and eq. 1.5.(19)]{Erd}

\[
{\int\limits_{0}^{\infty}  {\left( {1 + b\mu ^{z}} \right)\mu ^{x}d\mu}}   =
b^{ - {\displaystyle\frac{{x + 1}}{{z}}}}{\displaystyle\frac{{\Gamma \left( {{\displaystyle\frac{{x + 1}}{{z}}}}
\right)\Gamma \left( {y - {\displaystyle\frac{{x + 1}}{{z}}}} \right)}}{{z\Gamma \left(
{y} \right)}}}
\quad
{\left[ {b > 0,\,\,z > 0,\,\,\,\,0 < {\frac{{x + 1}}{{z}}} < y\,} \right]},
\]
and
\[
{\int\limits_{0}^{\pi / 2} {\sin ^{2x - 1}\mu \cos ^{2y - 1}\mu d\mu}}   =
{\frac{{\Gamma \left( {x} \right)\Gamma \left( {y} \right)}}{{2\Gamma \left(
{x + y} \right)}}}\,\,\,\,\,\,{\left[ {x > 0,\,\,y > 0} \right]}
\]
it is easy to obtain the equality (\ref{eq17}).

We note that a particular case (i.e. $t=0$) of the formula (\ref{eq17}) is found in the well-known  handbook \cite[p.635, eq. 4.638(3)]{grad}.

\section{Statement of the problem and an uniqueness theorem }

Consider equation (\ref{eq3}) in the infinite domain $\Omega $
defined in (\ref{eq4}).

We introduce the following notation:
\[
x:=\left(x_1,...,x_m\right)\in\mathbb{R}_m;\,\,\,R^{2}: = \sum_{i=1}^mx_{i}^{2};\,\,\, dx:=\prod_{i=1}^mdx_i;\,\,\, x^{(2\alpha)}:=\prod_{j=1}^nx_j^{2\alpha_j};
\]
\[
\tilde {x}_{k} : =\left( {x_{1} ,...,x_{k - 1} ,x_{k + 1}
,...,x_{m}}  \right) \in  \mathbb{R}_{m - 1};\,\,\, d\tilde {x}_{k} : =\frac{dx}{dx_k};\,\,\, \tilde{x}_{k}^{(2\alpha)}:=\frac{{x}^{(2\alpha)}}{x_k^{2\alpha_k}};
\]
\[
{x}_{k}^0 : =\left( {x_{1} ,...,x_{k - 1}, 0, x_{k + 1}
,...,x_{m}}  \right) \in  \mathbb{R}_{m};
\]

\[
S_{k} = \left\{ x:x_{1} > 0,...,x_{k - 1} > 0,\,x_{k} =
0, \,x_{k + 1} > 0,...,x_{n} > 0,  \right. \]
\[ \left.  - \infty < x_{n + 1} < +
\infty ,..., - \infty < x_{m} < + \infty  \right\},\,\,m\geq 2, \,\, 1\leq k \leq n\leq m.
\]

\textbf{The Neumann problem.} Find a regular solution  $u\left( {x} \right)$
of equation  (\ref{eq3}) from the class $C^1\left(\overline{\Omega}\right)\cap C^2(\Omega)$, satisfying the conditions:

\begin{equation}
\label{eq18}
{\left. {\left( {x_{k}^{2\alpha _{k}}  {\frac{{\partial u}}{{\partial x_{k}
}}}} \right)} \right|}_{x_{k} = 0} = \nu _{k} \left( {\tilde {x}_{k}}
\right),
\,\,
\tilde {x}_{k} \in S_{k} ,
\end{equation}

\begin{equation}
\label{eq19}
{\mathop {\lim} \limits_{R \to \infty}}  u\left( {x} \right) = 0,
\end{equation}
\noindent
where
 $\nu _{k} \left( {\tilde {x}_{k}}  \right)$  are given continuous functions.
The functions $\nu_k\left(\tilde {x}_{k}\right)$ can also turn to infinity of order less than $1-2\alpha_k$ and for sufficiently large values of $R$ the inequaliities are valid

\begin{equation}
\label{eq21}
{\left| {\nu _{k} \left( {\tilde {x}_{k}}  \right)} \right|} \le
{\frac{{c_{k}}} {{\left( 1+{ x_{1}^{2} + ... + x_{k - 1}^{2} + x_{k +
1}^{2} + ... + x_{m}^{2}}  \right)^{{\displaystyle{{\left(1 - 2\alpha _{k} + \varepsilon
_{k}\right)}}/{{2}}}}}}},
\end{equation}
where $c_{k} = const > 0$, $0 < 2\alpha _{k} < 1,$ and  $\varepsilon _{k} $
are small enough positive numbers $\left(k=\overline{1,n}\right)$.

\textbf{Theorem 1.}
 The Neumann problem has not more than one solution.

\textbf{Proof}. Suppose by contradiction. Let there be two solution $u_1$ and $u_2$ of the Neumann problem. We will denote via $u=u_1-u_2$. Then it is clear that the function $u$ satisfies the equation (\ref{eq3}) and homogeneous boundary conditions (\ref{eq18}) and condition (\ref{eq19}).

By ${D_R}$ we denote a bounded domain with a boundary $\partial D_R= \bigcup_{k=1}^nS_{Rk}$, where

\[
S_{Rk} = \left\{ x:0<x_{1}<R,...,0<x_{k - 1}<R,\,x_{k} =
0, \,0<x_{k + 1}<R,...,0<x_{n}<R,  \right. \]
\[ \left.  -R < x_{n + 1} < R ,..., -R < x_{m} <R  \right\}
\]
\[
\sigma _{R} : = \left\{ x:x_{1}^{2} + ... + x_{m}^{2} = R^{2},\,\,\,x_{1}
> 0,\,...,\, x_n>0\, - R < x_{n + 1} < +
R ,..., - R < x_{m} < R  \right\}.
\]

Choosing $R$ big enough, we integrate equation (\ref{eq3}) over the domain $D_R$, having previously multiplied it by the function $u(x)$, we obtain

\begin{equation}\label{eq21000}
\int_{D_R}x^{(2\alpha)}u\sum_{k=1}^m  \frac{\partial^2u}{\partial x_k^2}dx=0.
\end{equation}
Taking into account in (\ref{eq21000}) the following equalities
\[
u \frac{\partial^2u}{\partial x_k^2}=\frac{\partial}{\partial x_k}\left(u \frac{\partial u }{\partial {x_k}}\right)- \left(\frac{\partial u}{\partial {x_k}}\right)^2,\,\, k=\overline{1,n},
\]
after applying the Gauss-Ostrogaradsky formula, we have

\begin{equation} \label{eq21002}
\int_{D_R}x^{(2\alpha)}\sum_{i=1}^m\left(\frac{\partial u}{\partial x_k}\right)^2dx=\int_{\sigma_{R}}x^{(2\alpha)}u \frac{\partial u}{\partial {\mathbb{N}}}dS,
\end{equation}
where
\[
\frac{\partial u}{\partial {\mathbb{N}}}=\sum_{k=1}^n \frac{\partial u}{\partial {x_k}} \cos\left({\mathbb{N}{}},\, x_k\right);  \,\, \cos\left({\mathbb{N}},\, x_k\right)dS=d\tilde {x}_{k},\,\, k=\overline{1,n},
\]
${{\mathbb{N}}}$ is outer normal to $\partial D_R$.

By virtue of condition (\ref{eq19}), at $R \to \infty$, taking into account that \[{\mathop {\lim} \limits_{R \to \infty}}\int_{\sigma_{R}}x^{(2\alpha)}u \frac{\partial u}{\partial {\mathbb{{N}}}}dS=0\]
from (\ref{eq21002}), we have

\begin{equation} \label{eq21003}
\int_{D_R}x^{(2\alpha)}\sum_{i=1}^m\left(\frac{\partial u}{\partial x_i}\right)^2dx=0.
\end{equation}
From (\ref{eq21003}) we obtain $\displaystyle\frac{\partial u}{\partial x_i}=0,$ \,$(i=\overline{1,m})$ which means $u=const$. From the condition (\ref{eq19}) follows that $u\equiv 0$. So, we have proved the uniqueness theorem
for the Neumann problem.

\section{Existence of a solution of the Neumann problem}

\bigskip

Consider the function

\begin{equation}
\label{eq23}
u\left( {\xi}  \right):= u\left( {\xi_1,...,\xi_m}  \right) = - {\sum\limits_{k = 1}^{n} {{\int_{S_{k}}  {\tilde
{x}_{k}^{\left( {2\alpha}  \right)} \nu _{k} \left( {\tilde {x}_{k}}
\right){q \left( {x_k^0,\xi}  \right)} dS_{k}}
}\,}} ,
\end{equation}
where  $\nu _{k} \left( {\tilde {x}_{k}}  \right)$  are functions, defined in
 (\ref{eq18}) and  $q\left( {x,\xi}\right)$ is fundamental solution of equation (\ref{eq3})
\cite{{ergsib},{erglob}}:

\begin{equation*}
\label{eq24}
q \left( {x,\xi}  \right) = \gamma\,r^{ - 2\beta
}F_{A}^{\left( {n} \right)} \left({\beta, \alpha_1,...,\alpha_n;\, 2\alpha_1,...,2\alpha_n;- {\frac{4x_1\xi_1}{{r^{2}}}},...,- {\frac{4x_n\xi_n}{{r^{2}}}}}
\right),
\end{equation*}
\noindent
\begin{equation}
\label{eq25}
\beta:=\displaystyle\frac{m-2}{2}+\alpha,\,\,
\gamma = 2^{2\beta - m}{\frac{{\Gamma \left(\beta
\right)}}{{\pi ^{m / 2}}}}{\prod\limits_{k = 1}^{n} {{ {{\frac{{\Gamma
\left( {\alpha _{k}}  \right)}}{{\Gamma \left( {2\alpha _{k}}  \right)}}}}
}}} ,\,\, r^{2} = {\sum\limits_{i = 1}^{m} {\left( {x_{i} - \xi _{i}}  \right)^{2}}}.
\end{equation}

Here

\[
{\int_{S_{k}}  {f\left( {x,\xi}  \right)dS_{k}}}  : = {\underbrace
{{\int\limits_{0}^{ + \infty}  {...{\int\limits_{0}^{ + \infty}  {}}} }} _{n
- 1}}{\underbrace {{\int\limits_{ - \infty} ^{ + \infty}  {...{\int\limits_{
- \infty} ^{ + \infty}  {}}} }} _{m - n}}f\left( {x,\xi}  \right)dx_{1}
...dx_{k - 1} dx_{k + 1} ...dx_{n} dx_{n + 1} ...dx_{m} .
\]

It is easy to check that the fundamental solution $ q \left ({x, \xi}
\right) $ has the properties

\[
{\left. {{\frac{{\partial q \left( {x,\xi}  \right)}}{{\partial x_{k}
}}}} \right|}_{x_{k} = 0} = 0, \,\,\,k = \overline {1,n}.
\]

Let us prove that the function (\ref {eq23}) is a solution to the Neumann problem. It is clear that the function (\ref {eq23}) in the arguments $\xi$ satisfies the equation (\ref{eq3}) (for details, see ). Let us show that the function (\ref{eq23}) also satisfies the conditions of the Neumann problem. Due to the fact that ${\left.{\sigma_{k}}\right|}_ {x_ {k} = 0} = 0$, the expression (\ref {eq23}) is converted to the form

\begin{equation}
\label{eq26}
u\left( {\xi}  \right) = {\sum\limits_{j = 1}^{n} {\,I_{j} \left( {\xi}
\right)}} ,
\end{equation}
where

\begin{equation}
\label{eq27}
I_{j} \left( {\xi}  \right) = - {\gamma} {\int_{S_{j}}  {\nu _{j}
\left( {\tilde {x}_{j}}  \right)}} \,\tilde {x}_{j}^{\left( {2\alpha}
\right)}{r}_{j}^{ - 2{\beta}}  F_{A}^{\left( {n - 1} \right)}
{\left( {{{{{{{\beta}}}} ,{{\rm{A}}}_j ;} {2{{\rm{A}}}_{j};}
}{\Phi}_{j}}\right)}dS_{j} ,\,
\end{equation}

\[
{\rm{A}}_j:= \left(\alpha_1,...,\alpha_{j-1}, \alpha_{j+1},...,\alpha_n\right),\,\,{\Phi}_j:=\left(\phi_{1,j},...,\phi_{j-1,j}, \phi_{j+1,j},...,\phi_{n,j}\right),
\]
\[
 \phi_{i,j}:=-\frac{4x_{i}\xi_{i}}{{r}_j^2}, \,\,i\neq j;\,\,\,{r}_{j}^{2}: = {\left. r^2\right|}_{x_j=0},\,\,  i,j=\overline{1,n}.
\]

Let us show that the function (\ref{eq27}) satisfies the conditions (\ref{eq18}) of the Neumann problem as
for $j=k$ and for $ j\ne k $.

Let $j=k$. Using the differentiation formula (\ref{differential}) and adjacent relation (\ref{adjacent}) for Lauricella
hypergeometric  function, we obtain

\begin{equation}
\label{eq30}
{\frac{{\partial I_{k}}} {{\partial \xi _{k}}} } = 2\beta { \gamma
} \xi _{k} {\int_{S_{k}}  {\nu _{k} \left( {\tilde {x}_{k}}  \right)}
}\,\,r_{k}^{ - 2{\beta} - 2} \tilde {x}_{k}^{\left( {2\alpha}
\right)} F_{A}^{\left( {n - 1} \right)} {\left( {{
 {\beta}+ 1, {\rm{A}}_k ;  2{\rm{A}}_k;
}  {\Phi}_{k}}  \right)}dS_{k} .
\end{equation}

On the right-hand side of the equality (\ref{eq30}), we make the change of variables

\[
x_{i} = \xi _{i} + \xi _{k} t_{i} ,\,\,i = \overline {1,m} ,\,\,\,i \ne k,
\]
then we have

\[
\xi _{k}^{2\alpha _{k}}  {\frac{{\partial I_{k}}} {{\partial \xi _{k}}} } =
2\beta{\gamma} \,{\int_{T_{k}}  {{\frac{{\nu _{k} \left(
{\xi _{1} + \xi _{k} t_{1} ,...,\xi _{k - 1} + \xi _{k} t_{k - 1} ,\xi _{k +
1} + \xi _{k} t_{k + 1} ,...,\xi _{m} + \xi _{k} t_{m}}  \right)}}{{\left(
{1 + t_{1}^{2} + ... + t_{k - 1}^{2} + t_{k + 1}^{2} + ... + t_{m}^{2}}
\right)^{\beta+1}}}}}} \,
\]

\[
 \times \prod\limits_{{s = 1,}{s \ne k}}^{n}\left(\frac{\xi_s+\xi_kt_s}{\xi_k}\right)^{2\alpha_s}F_{A}^{\left( {n - 1} \right)} \left[\begin{array}{*{20}c}
 {\beta}+ 1, \alpha_1,...,\alpha_{k-1}, \alpha_{k+1},...,\alpha_n; \hfill\\
 2\alpha_1,...,2\alpha_{k-1}, 2\alpha_{k+1},...,2\alpha_n;\hfill\\
 \end{array}1-\frac{z_1\left(\tilde{t}_k,\xi_1,\xi_k\right)}{\xi_k^2},...\,\,\,\,\,\,\,\,\,\,\, \right.
\]

\begin{equation}
\label{eq35}
 \left.\begin{array}{*{20}c} {} \hfill\\
 {} \hfill\\ \end{array}\,\,\,\,\,\,\,\,\,\,\,\,\,\,\,\,\,\,\,\,\,\,\,\,\,\,\,\,\,\,\,\,\,...,1-\frac{z_{k-1}\left(\tilde{t}_k,\xi_{k-1},\xi_k\right)}{\xi_k^2}, 1-\frac{z_{k+1}\left(\tilde{t}_k,\xi_{k+1},\xi_k\right)}{\xi_k^2},..., 1-\frac{z_n\left(\tilde{t}_k,\xi_n, \xi_k\right)}{\xi_k^2}
  \right] dT_{k} ,
\end{equation}
where
\[
z_{i}\left(\tilde{t}_k,\xi_{i},\xi_k\right):= \frac{4\xi_i\left(\xi_i+\xi_kt_i\right)+\left(1 + t_{1}^{2} + ... + t_{k - 1}^{2} + t_{k + 1}^{2} + ... + t_{m}^{2}\right)\xi_k^2}{1 + t_{1}^{2} + ... + t_{k - 1}^{2} + t_{k + 1}^{2} + ... + t_{m}^{2}},
\]

\[
\,\tilde {t}_{k} : = \left( {t_{1} ,...,t_{k - 1} ,t_{k + 1} ,...,t_{n}
,...,t_{m}}  \right),
\]

\[
{\int_{T_{k}}  {...dT_{k}}}  : = {\underbrace {{\int\limits_{ - {\frac{{\xi
_{1}}} {{\xi _{k}}} }}^{ + \infty}  {...{\int\limits_{ - {\frac{{\xi _{k -
1}}} {{\xi _{k}}} }}^{ + \infty}  {\,{\int\limits_{ - {\frac{{\xi _{k + 1}
}}{{\xi _{k}}} }}^{ + \infty}  {...}}} } {\int\limits_{ - {\frac{{\xi _{n}
}}{{\xi _{k}}} }}^{ + \infty}  {}}} }} _{n - 1}}{\underbrace {{\int\limits_{
- \infty} ^{ + \infty}  {...{\int\limits_{ - \infty} ^{ + \infty}  {}}}
}}_{m - n}}...dt_{1} ...dt_{k - 1} dt_{k + 1} ...dt_{n} dt_{n + 1} ...dt_{m}
.
\]

In the both sides of the equality (\ref{eq35}), we pass to the limit as $ \xi_{k} \to 0$. Taking into account formulae (\ref{eq16}) and (\ref{eq17}),  \,Legendre's duplication formula  (\ref{legendre})
and the expression for the coefficient $ \gamma$  by (\ref{eq25}), we have

\begin{equation}
\label{eq39}
{\mathop {\lim} \limits_{\xi _{k} \to 0}} \xi _{k}^{2\alpha _{k}}
{\frac{{\partial I_{k}}} {{\partial \xi _{k}}} }\, = \nu _{k} \left( {\tilde
{\xi} _{k}}  \right).
\end{equation}

Let now $l \ne k(l,k = \overline {1,n} )$.  For definiteness, we put $ l <k $.

 Repeating the reasoning for obtaining the formula (\ref{eq30}), we have

\[
{\frac{{\partial I_{k}}} {{\partial \xi _{l}}} } = 2\beta { \gamma
} {\int_{S_{k}}  {\nu _{k} \left( {\tilde {x}_{k}}  \right)}} \,\,\left(
{x_{l} - \xi _{l}}  \right)r_{k}^{ - 2{\beta - 2}} \tilde
{x}_{k}^{\left( {2\alpha}  \right)} F_{A}^{\left( {n - 1} \right)} {\left[
{\begin{array}{*{20}c}
 {\beta + 1,{\alpha} _{1} ,...,{
\alpha} _{k - 1} ,{\alpha} _{k + 1} ,...,{\alpha} _{n} ;} \hfill \\
 {{2}{\alpha} _{1} ,...,2{\alpha} _{k - 1} ,2{\alpha} _{k +
1} ,...,2{\alpha} _{n} ;} \hfill \\
\end{array} \tilde {\sigma} _{k}}  \right]}dS_{k}
\]

\[
 - 2\beta {\gamma} {\int_{S_{k}}  {\nu _{k} \left( {\tilde
{x}_{k}}  \right)}} \,\,r_{k}^{ - 2{\beta - 2}} x_{l} \tilde
{x}_{k}^{\left( {2\alpha}  \right)} F_{A}^{\left( {n - 1} \right)} \left[
{\begin{array}{*{20}c}
 {\beta + 1,{\alpha} _{1} ,...,\alpha _{l}+ {
1},...,{\alpha} _{k - 1} ,{\alpha} _{k
+ 1} ,...,{\alpha} _{n} ;} \hfill \\
 {{2}{\alpha} _{1} ,...,2\alpha _{l} + 1,...,2{\alpha}_{k - 1}
,2{\alpha} _{k + 1} ,...,2{\alpha} _{n} ;} \hfill \\
\end{array} \tilde {\sigma} _{k}}  \right]dS_{k} .
\]

Now using the obvious equality

\[
{\left. {\left( {x_{l} - \xi _{l}}  \right)F_{A}^{\left( {n - 1} \right)}
{\left[ {{\begin{array}{*{20}c}
 {{1}{+}\beta,{\alpha}_{1} ,...,\alpha_l,...,{\alpha}_{k
- 1} ,{\alpha}_{k + 1} ,...,{\alpha}_{n};} \hfill \\
 {{2}{\alpha}_{1} ,...,2\alpha_l, ..., 2{\alpha}_{k - 1} ,2{\alpha}_{k +
1} ,...,2{\alpha}_{n} ;} \hfill \\
\end{array}} \Phi_k}  \right]}} \right|}_{\xi _{l} = 0} =
\]

\[
{\left. { = x_{l} F_{A}^{\left( {n - 1} \right)} {\left[
{{\begin{array}{*{20}c}
 {{1}{+}\beta, {\alpha}_{1},...,1 + \alpha_{l}
,...,{\alpha}_{k - 1} ,{\alpha} _{k +
1} ,...,{\alpha} _{n} ;} \hfill \\
 {{2}{\alpha}_{1} ,...,1 + 2\alpha_{l} ,...,2{\alpha}_{k - 1}
,2{\alpha}_{k + 1} ,...,2{\alpha}_{n} ;} \hfill \\
\end{array}} \Phi_k}  \right]}} \right|}_{\xi _{l} = 0}
\]

\noindent
it is easy to prove that

\begin{equation} \label{lll}
{\mathop {\lim} \limits_{\xi _{l} \to 0}} \xi _{l}^{2\alpha _{l}}
{\frac{{\partial I_{k}}} {{\partial \xi _{l}}} }\, = 0,\,\,\,l < k.
\end{equation}

In this way,

\begin{equation}
\label{eq40}
{\mathop {\lim} \limits_{\xi _{l} \to 0}} \xi _{l}^{2\alpha _{l}}
{\frac{{\partial I_{k}}} {{\partial \xi _{l}}} }\, = 0,\,\,\,l > k.
\end{equation}

Therefore, by virtue of equalities (\ref{eq39}), (\ref{lll}) and (\ref{eq40}), we conclude that the function (\ref{eq23}) satisfies the conditions (\ref{eq18}) of the Neumann problem.

Next, we show that if the given functions $ \nu_{k} \left({\tilde{x}_{k}} \right) $ for sufficiently large values of the arguments satisfy the inequalities (\ref{eq21}), then the solution (\ref{eq23}) of the Neumann problem also satisfies the condition (\ref{eq19}).

 Indeed, let the inequalities (\ref{eq21}) hold, then in the equalities (\ref{eq27}) we make the following changes of variables

\[
y_{i} = {\frac{{1}}{{R_{0}}} }x_{i} ,
\quad
\eta _{i} = {\frac{{1}}{{R_{0}}} }\xi _{i} ,
\quad
i = \overline {1,m} ,
\]

\noindent
where $R_{0}^{2} : = \xi _{1}^{2} + ... + \xi _{m}^{2} .$

Then, by virtue of (\ref{eq27}) we have the following inequalities for $ R_ {0} \to \infty $

\begin{equation}
\label{eq41}
{\mathop {\lim} \limits_{R_{0} \to \infty}}  {\left| {I_{k} \left( {\xi}
\right)} \right|} \le {\frac{{{\rm 2}^{m - n}{\rm \gamma} c_{k}
}}{{R_{0}^{\varepsilon _{k}}} } }{\underbrace {{\int\limits_{0}^{ + \infty}
{...{\int\limits_{0}^{ + \infty}  {}}} }} _{m - 1}}{\frac{{dy_{1} ...dy_{k -
1} dy_{k + 1} ...dy_{m}}} {{ {Y}^{1 - 2\alpha _{k} + \varepsilon _{k}
}\left( {1 + {Y}^{2}} \right)^{\beta}} }}\,,
\end{equation}

\noindent
where ${Y}^{2}: = y_{1}^{2} + ... + y_{k - 1}^{2} + y_{k + 1}^{2} + ...
+ y_{m}^{2} .$

Now we will show that $(m - 1)-$ dimensional integrals occurring in inequalities (\ref{eq41}),
are limited.

Indeed, by virtue of the value of the integral (\ref{eq17}), the inequalities (\ref{eq41}) imply

\begin{equation}
\label{eq42}
{\mathop {\lim} \limits_{R_{0} \to \infty}}  {\left| {I_{k} \left( {\xi}
\right)} \right|} \le {\frac{{\tilde {c}_{k}}} {{R_{0}^{\varepsilon _{k}}
}}}\,,
\end{equation}

\noindent
where $\tilde {c}_{k}  $ are constants. Due to the inequalities (\ref{eq42}) from the expression
(\ref{eq26}), we finally get

\[
{\left| {I_{k} \left( {\xi}  \right)} \right|} \le
{\frac{{c}}{{R_{0}^{\varepsilon}} } }\,,\,
c = const,
\,\varepsilon = {\mathop {\min} \limits_{1 \le k \le n}} {\left\{
{\varepsilon _{k}}  \right\}}.
\]

The last inequality shows that the solution (\ref{eq23}) vanishes for $ R_{0} \to \infty $. Therefore, condition (\ref{eq19}) of the Neumann problem is satisfied. Thus, the solution (\ref{eq23}) satisfies all the conditions of the Neumann problem.

\end{document}